\documentclass[final]{siamart220329}

\usepackage{amsmath, amssymb}
\usepackage{mathtools, commath}
\usepackage{stmaryrd}
\usepackage{booktabs}
\AtBeginEnvironment{tabular}{\footnotesize}
\newcommand{\Reyuls}{\mathbb{R}}
\newcommand{\bo}{{\partial\Omega}}
\newcommand{\xx}{{\mathbf x}}
\newcommand{\vv}{{\mathbf v}}
\newcommand{\uu}{{\mathbf u}}
\newcommand{\ww}{{\mathbf w}}
\newcommand{\nn}{{\mathbf n}}
\newcommand{\bfz}{{\mathbf 0}}
\newcommand{\gbc}{{\mathbf g}}
\newcommand{\sdiv}{{{\nabla\cdot} \,}}
\newcommand{\du}{{\cal D}(\uu)}
\newcommand{\dv}{{\cal D}(\vv)}

\newcommand{\half}{{\textstyle{\frac{1}{2}}}}
\newcommand{\threehalfs}{{\textstyle{\frac{3}{2}}}}

\newcommand{\hdiv}{H(\text{div})}
\newcommand{\hcurldiv}{H(\text{curl div})}
\newcommand{\bog}{\Gamma}
\newcommand{\famlyT}{{\mathcal T}_h}
\newcommand{\famlyF}{{\mathcal F}_h}
\newcommand{\SVfem}{\texttt{SV}}
\newcommand{\MCSfem}{\texttt{MCS}}
\newcommand{\gdTHfem}{\texttt{gdTH}}
\newcommand{\THfem}{\texttt{TH}}

\newcommand{\jump}[1]{\llbracket#1\rrbracket}

\newsiamremark{remark}{Remark}

\usepackage[textsize=footnotesize,textwidth=5cm]{todonotes}
\makeatletter
\@mparswitchfalse%
\makeatother
\normalmarginpar

\title{Reliable chaotic transition in incompressible fluid simulations\thanks{%
Part of this research was conducted using computational resources and services
at the Center for Computation and Visualization, Brown University.
\funding{This material is based upon work supported by the National Science
Foundation under Grant No.~DMS-1929284 while the authors were in residence at
the Institute for Computational and Experimental Research in Mathematics in
Providence, RI, during the Numerical PDEs: Analysis, Algorithms, and Data
Challenges program.}}
}

\author{Henry von Wahl%
\thanks{Institute for Computational and Experimental Research in Mathematics,
Brown University, 121 South Main Street, Box E, Providence, RI 02903, USA, and
Fakultät für Mathematik und Informatik, Friedrich-Schiller-Universität,
Ernst-Abber-Platz 2, 07743 Jena, Germany, \email{henry.vonwahl@uni-jena.de}}
\and
L. Ridgway Scott%
\thanks{Department of Mathematics, University of Chicago, Chicago,
Illinois 60637, USA, \email{ridg@cs.uchicago.edu}}
}

\begin{document}

\maketitle

\begin{abstract}
We consider a test problem for Navier--Stokes solvers based on the flow around
a cylinder that exhibits chaotic behavior, to examine the behavior of various
numerical methods. We choose a range of Reynolds numbers for which the flow is
time-dependent but can be characterized as essentially two-dimensional. The
problem requires accurate resolution of chaotic dynamics over a long time
interval. It also requires the use of a relatively large computational domain,
part of which is curved. We review the performance of
different finite element methods for the proposed range of Reynolds numbers.
These tests indicate that some of the most established methods do not capture
the correct behavior. The key requirements identified are pressure-robustness
of the method, high resolution, and appropriate numerical dissipation when
the smallest scales are under-resolved.
\end{abstract}



\section{Introduction}

We consider a classic fluids problem, the chaotic flow around a cylinder, as a
challenge for numerical schemes. As reviewed in \cite{lrsBIBki}, there is ample
experimental data with which to compare simulations. To make the challenge more
tractable, we focus on 2D simulations and a range of Reynolds numbers.
Illustrating the difficulty posed by this problem,
we review some published examples where erroneous results were obtained. We
also show that some commonly used techniques fail to give acceptable results
at Reynolds numbers of interest.

Flow around a cylinder is not only a fundamental problem in fluid dynamics but
also a practical problem of interest in energy generation. For example,
so-called bladeless turbines \cite{bahadur2022dynamic,el2018vortex,francis2021design}
have been proposed as a viable energy generation method. The design of such
systems requires accurate and efficient simulations of the flow around
cylinders. Another area of recent development in which such simulations are
necessary is new design spaces for electric aircafts, such as ground-effect
transportation~\cite{dBre18} and air taxis~\cite{RS20}.

As computational simulation has emerged as a
new form of experimentation in physics, it is essential to have some benchmarks
that give guidance, as we present here. For many flow problems, there is no
dispute regarding the Navier--Stokes equations as a suitable model, as we do
here. However, this requires reliable methods to decide fundamental physics
questions, such as whether the resulting flow is chaotic at a given Reynolds number.

A critical factor in obtaining reliable results at high Reynolds numbers and
avoiding erroneous chaos, is for a numerical method is to
have an appropriate mechanism to deal with the energy from under resolved
scales~\cite{Bur07}. This usually takes the form of numerical dissipation and
can come from different stabilization terms in the discrete method. However,
it is also critical for the method not to be overly dissipative, which would
again yield erroneous results. In this work, we cover a number of schemes with
different dissipative mechanisms providing significant data clarifying this.

Flow around a cylinder has been proposed as a test problem
before~\cite{schafer1996benchmark,ref:refvalcyliftdragVolkerJohn}. What is
different here is that we focus on a range of Reynolds numbers for which the
flow appears to be chaotic \cite{lrsBIBki}. Thus, appropriate metrics must be
used to compare different numerical methods. In particular, we focus on
statistics predicting the periodicity of the flow, i.e. the Strouhal
period~\cite{Str78}, rather than just minimal and maximal drag and lift
coefficient values. Moreover, the time interval for the simulation is quite
long, so this simulation evaluates the ability of different discretizations
with regard to accumulation of numerical error and challenges the
implementations with regard to efficiency.

For a straightforward analysis of the resulting flow,
we focus on simple metrics so that different methods can be compared easily.
The Strouhal number or period of the time-dependent flow is such a quantity.
This can be computed by analyzing the phase-diagram of the drag and lift
coefficients, which are computed easily. It is vital to compute the period by
considering the drag and lift simultaneously, as both have local maxima and
minima making it difficult to determine the period in cases where the flow is
not perfectly periodic. In \cite{lrsBIBki}, a simple approach was
used to estimate the periodicity of chaotic flow. We present a slightly
improved method compared to \cite{lrsBIBki} for evaluating the Strouhal period
that provides an indication of the transition to chaos. Additional metrics,
such as the Lyapunov exponent and the fractal dimension of the attractor,
provide a finer level of detail regarding the flow dynamics \cite{lrsBIBki}.

Evaluating computational methods for fluid flows is already viable in two
spatial dimensions and sufficiently challenging to differentiate between
them. We, therefore, focus on two dimensional flows. When the Reynolds number
is more than $10^4$, the flow becomes three-dimensional
\cite{ref:FazleHayakawatreedee}. Thus, we restrict ourselves to Reynolds numbers
in the range $[10^2,10^4]$ for simplicity. Three-dimensional effects of a
different kind are seen \cite{ref:obliqueparallelmodesWilliamson} at lower
Reynolds numbers due to the finite length of cylinders used in experiments.
However, these are potentially like so-called blocking effects in
two-dimensional flow \cite{ref:blockagecylinderatio}, which become less
significant at higher Reynolds numbers.

In this paper, we consider a variety of well-studied finite element methods
with a solid theoretical foundation. These are low and high-order Taylor--Hood,
with and without grad-div stabilization, Scott--Vogelius finite elements and a
mass-conserving mixed stress (MSC) based formulation using $\hdiv$-conforming
velocities. The latter we consider both with an upwind and central flux in
the convective term. It is rare for different numerical methods to be compared
head-to-head on a challenging physical problem of interest.
We hope that this paper will empasize the value of doing so.
Although we have studied several commonly used methods, and we compare results
with other published simulation methods, there are many more methods that have
been proposed for simulations of the type we do here.
We hope our benchmark problem will attract others to determine to what extent
their methods compare with the results presented here. Finally, we compare
their performance to a method used in \cite{lrsBIBki}, and we find that the
simulations in \cite{lrsBIBki} were under-resolved for higher Reynolds numbers.

\begin{figure}
\centering
\includegraphics{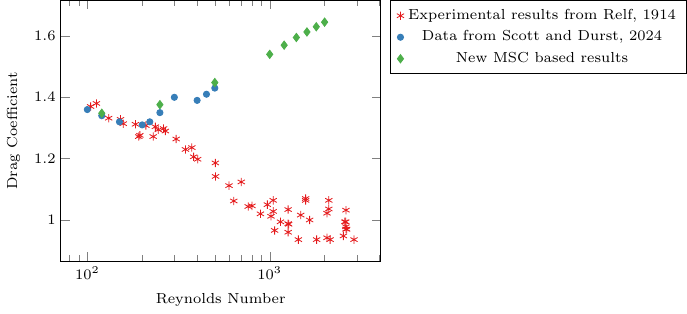}
\caption{Relf experimental data \cite{ref:errelflowRex} compared to the data computed
in \cite{lrsBIBki} with the Scot--Vogelius finite elements scheme and the new data presented here.}
\label{fig:relfdata}
\end{figure}

Our main conclusion is that
high-order methods are necessary to obtain reliable results at higher Reynolds
numbers. Surprisingly, the commonly used, lowest-order Taylor--Hood method is
not accurate enough to be used even at the lowest Reynolds numbers considered
here. A secondary conclusion is that it appears that methods that result in
exactly divergence-free velocity  solutions have a significant advantage at
higher Reynolds numbers, even though one that we investigate is non-conforming.
Furthermore, we also conclude that a viable method at high Reynolds numbers
needs to contain the correct amount of numerical dissipation to capture the
flow dynamics on coarser meshes. Finally, one of the main conclusions from
\cite{lrsBIBki} was that the simulation drag values differ substantially from
experiments, as shown in \Cref{fig:relfdata}. We confirm this conclusion here,
suggesting that cylinder vibration likely plays a significant role in practice
\cite{ref:annurevfluidvortexvibrations}.

The remainder of this paper is structured as follows. In \Cref{sec.equations}, we discuss the equations and the quantities under consideration. In \Cref{sec.setup}, we then present the details of the computational setup for our computations. \Cref{sec.disctete-methods} covers the details of the numerical methods under consideration. In \Cref{sec.numex}, we present and discuss the numerical results achieved. In particular, we start at a low Reynolds number to illustrate how established methods struggle to achieve the correct dynamics. We then cover a range of Reynolds numbers up to $Re=8000$. Finally, we give some concluding remarks in \Cref{sec.conclusions}.

\section{Setting the problem and model equations}
\label{sec.equations}

Suppose that $(\uu,p)$ is a solution of the time-dependent Navier--Stokes
equations in a domain $\Omega\subset\Reyuls^d$ containing an obstacle with
boundary $\Gamma\subset\bo$. This fulfills the equations
\begin{equation} \label{eqn:firstnavst}
\begin{split}
\partial_t\uu-\nu\Delta \uu + \uu\cdot\nabla\uu+\nabla p &= \bfz\quad\hbox{in}\;\Omega,\\
\sdiv\uu &=0\quad\hbox{in}\;\Omega,
\end{split}
\end{equation}
with the kinematic viscosity $\nu$, and together with boundary conditions
\begin{equation} \label{eqn:bceesnavst}
\uu=\gbc\;\hbox{on}\;\partial\Omega\backslash\Gamma,\quad
\uu=\bfz\;\hbox{on}\;\Gamma.
\end{equation}
For the well-posedness of these equations, see \cite{giraultraviart}.

\subsection{Weak formulation of the Navier--Stokes equations}
While this test problem is not limited to finite-element
simulations, we expect many will be done with finite elements. To this end,
we require a variational formulation of \eqref{eqn:firstnavst}.
The Navier--Stokes equations can be written in a weak (or variational)
form as follows:

\smallskip
Find $\vv \in H^1_g(\Omega)$ and $p \in L^2_0(\Omega)$ such that
\begin{equation}\label{eqn:varform}
\begin{split}
(\partial_t\uu,\ww)_{L^2(\Omega)}+ a(\uu, \ww) + b(\ww, p) +c(\uu,\uu,\ww)&= 0
\quad \forall \ww \in (H^1_0(\Omega))^d, \\
b(\uu, q) &= 0 \quad \forall q \in L^2_0(\Omega).
\end{split}
\end{equation}
The space $H^1_g(\Omega)$ is the space of vector-valued $H^1$ functions with
trace $g$ on the (outer) boundary, $H^1_0(\Omega)$ is the space of $H^1$
functions with trace $0$ on the boundary and $L^2_0(\Omega)$ is the space of
$L^2$ functions with mean zero. The bilinear forms $a(\cdot, \cdot)$ and
$b(\cdot, \cdot)$ are defined as
\begin{equation*}
a(\vv, \ww) \coloneqq \int_\Omega \frac{\nu}{2} \mathcal{D}(\vv):\mathcal{D}(\ww)\dif\xx,
\quad\text{and}\quad
b(\vv, q) \coloneqq -\int_\Omega q (\nabla \cdot \vv) \dif\xx,
\end{equation*}
respectively, where $\dv=\nabla\vv+\nabla\vv^t$ and the colon (:) indicates the
Frobenious inner-product of matrices. The convective (nonlinear) term takes various forms.
By default, we take it to be
\begin{equation*}
   c(\uu,\vv,\ww)=c_{conv}(\uu, \vv, \ww) \coloneqq \int_\Omega (\uu\cdot\nabla \vv)\cdot \ww  \dif\xx .
\end{equation*}
For divergence-free $H^1$ functions, the various forms of the convective term
are equivalent. Hence, we may take the above so-called convective form
of the trilinear form for exactly divergence-free methods. However, for
non-divergence-free methods, other forms are beneficial and for
non-$H^1$-conforming methods are necessary.

\subsection{Drag and lift evaluation}

The flow of fluid around an obstacle generates a force called drag, which is a
fundamental concept in fluid dynamics \cite{lamb1993hydrodynamics}.
It plays a critical role in determining the behavior of objects in flight and
has been studied since the time of d'Alembert \cite{lrsBIBjn}.
Drag is composed of two components: pressure drag $\beta_{p}$ and
viscous drag $\beta_{v}$, which can be calculated by evaluating the
following functions using $\vv=(1,0)$:
\begin{equation*}
   \beta_v(\vv) =\int_{\bog} \big((\nu\du)\vv\big)\cdot\nn\dif s,\quad
   \beta_p(\vv) =\int_{\bog}  -p\vv\cdot\nn\dif s.
\end{equation*}
The full drag $\beta$ is defined by $\beta=\beta_{v}+\beta_{p}$.
Similarly, lift is computed using $\vv=(0,1)$.

\subsubsection{Alternate drag evaluation}

Another way \cite{ref:refvalcyliftdragVolkerJohn} to evaluate $\beta$ is
to test the weak formulation \eqref{eqn:varform} with a non-conforming test
function. This results in the functional
\begin{equation*}
   \omega(\vv)=
     \int_\Omega \partial_t\uu\,\vv\dif\xx
     +\int_\Omega \frac{\nu}{2}\du:\dv +(\uu\cdot\nabla\uu)\cdot\vv-p\sdiv\vv\dif\xx,
\end{equation*}
Then $\omega(\vv)=\beta(\vv)$ for all $\vv\in H^1(\Omega)^d$
\cite{ref:refvalcyliftdragVolkerJohn,lrsBIBkj}.
This approach may be traced back to \cite{babuvska1984post}. It is also known
that, for the case of strongly imposed Dirichlet boundary conditions, this
approach of testing the residual with a non-conforming test-function doubles
the rate of convergence for the drag and lift, see, for
example~\cite{braack2006solutions} for the proof in the steady case.

Having two ways to compute drag and lift provides a valuable internal
check of the accuracy of the simulation.

\section{Problem Set-up}
\label{sec.setup}

For our simulations, we consider the domain
\begin{equation*}
\Omega=\{(x,y):-30<x<300,\,|y|<30,\;x^2+y^2>1\}.
\end{equation*}
The boundary condition is set as $\gbc = (1, 0)^T$ and the cylinder diameter
is the reference length, so the Reynolds number is given by $Re = 2 / \nu$.
We consider the time interval $[0, 500]$.

\subsection{Challenge computation}
The data in \Cref{fig:relfdata} suggests a significant disagreement
between computation and experiment.
This could be for a variety of reasons.
The most obvious is some flaw in the computational scheme.
Thus, we propose this problem as a computational challenge to see
what other methods predict.

Another reason for the disagreement could be that the model is wrong.
Some people who have published experiments (see \cite{lrsBIBki})
have raised the issue that the cylinders used could vibrate.
In particular, the Relf data in \Cref{fig:relfdata} is based on flows past
thin wires, used in musical instruments for their vibrational qualities.
Thus, one response to the challenge is to allow the cylinder to vibrate
back and forth in response to the oscillating lift and drag.
The question would be to evaluate how oscillations of the cylinder change
parameters describing the flow, such as drag, lift, and Strouhal number,
representing the frequency of the generated vortices.
One study with forced vibrations \cite[Figure 6.9, page 179]{ref:Meneghini1994Thesis}
suggests that the change could be substantial.

\section{Discretization methods}
\label{sec.disctete-methods}

The simulations in this paper are done using finite element methods (FEM). We,
therefore, begin by presenting some details of the finite element approach.

To discretize the weak formulation \eqref{eqn:varform} with finite element
methods, we take a mesh of the domain $\Omega$ with characteristic length
$h$, denoted as $\famlyT$. Finite element spaces then discretize
\eqref{eqn:varform} by approximating the spaces $H^1(\Omega)$ and
$L^2(\Omega)$ with spaces of piecewise polynomials on each element of the
mesh, with varying degree of discontinuity of these polynomials over
element edges.

We utilize many well-established simulation methods for which extensive
numerical analysis is available, establishing not only convergence but
also detailed simulation properties.

\subsection{H(div)-conforming FEM}
\label{sec:disc.subsec:hdiv}

We consider an $\hdiv$-conforming finite element method based on
\cite{gopalakrishnan2020weak,gopalakrishnan2020mass}. This discretization
is based on a mass-conserving mixed formulation with symmetric stresses,
known as the MCS formulation. This approach discretizes the velocity in $\hdiv$
using Brezzi--Douglas--Marini (BDM) elements \cite{lrsBIBih}. That is, the
normal component of the velocity is continuous across element boundaries, while
the tangential component is discontinuous. The pressure is discretized using
discontinuous elements. This method is exactly mass conserving and
pressure-robust.
From a numerical point of view, pressure-robust methods have several advantages.
First, the velocity remains unchanged by irrotational changes in any external
forcing term, a property that most methods do
not preserve on the discrete level \cite{Link2_helmholtz_2014}.
More importantly, pressure-robust methods have the property that the velocity
error is independent of the pressure error~\cite{JLM17}.
This contrasts classical Galerkin discretizations of the
(Navier)-Stokes equations, where the velocity error depends on
$\frac{1}{\nu}$ times the pressure error. From a physical point of view,
exactly divergence-free methods have the advantage that, in addition to mass,
kinetic energy, linear momentum, and angular momentum are conserved on the
discrete level~\cite{CHOR17,SL17}.

This formulation introduces an additional variable
$\sigma = \nu\nabla\uu$ which is discretised in the function space $\hcurldiv$.
Thus the energy form is modified to
\begin{equation*}
a'(\sigma, \tau) \coloneqq \frac1\nu \int_\Omega \sigma:\tau \dif\xx,
\end{equation*}
and a new form is defined by
\begin{equation*}
b'(\tau,\uu) \coloneqq \sum_{T\in\famlyT} \int_{T} \sdiv(\tau)\cdot\uu \dif\xx
+\sum _{F\in\famlyF} \int_{F} \jump{\tau_{nn}}\uu\cdot\nn\dif s,
\end{equation*}
where $\famlyT$ is the triangular mesh, and $\famlyF$ is the set of element edges.
The expression $\jump{\tau_{nn}}$ denotes the jump in the normal-normal component
of the tensor $\tau$ across the edge $e$, and $\nn$ is the normal to the edge.
Note that the orientation of the normal is not material. As the velocity is
$\hdiv{}$, we have to modify the convective term to account for the discontinuity
of the velocity in the tangential direction, and there are multiple choices available.
For the most part, we consider the upwind flux for the convective term,
\begin{equation}
   c_{upw}(\uu_h,\uu_h,\ww_h)
      \coloneqq \sum_{T\in\famlyT} \int_T -(\uu_h\cdot\nabla \ww_h)\cdot \uu_h \dif\xx
      + \int_{\partial T} \uu_h\cdot\nn\, \hat{\vv}_h\cdot\ww_h \dif s,
\end{equation}
where $\hat{\vv}_h$ is the upwind flux, i.e., $\hat{\vv}_h =\uu_h$ if
$\uu_h(x)\cdot\nn(x)\geq 0$ and otherwise $\hat{\vv}_h =\uu_h^\text{oth}$,
where $\uu_h^\text{oth}$ is the value of $\uu_h$ from the neighboring element.
An alternative choice is a central flux. In this case the convective term can
be written as
\begin{equation}
   c_{cf}(\uu_h,\uu_h,\ww_h)
    \coloneqq \sum_{T\in\famlyT} \int_T (\uu_h\cdot\nabla \uu_h)\cdot \ww_h \dif\xx
      - \int_{\partial T} \uu_h\cdot\nn (\hat{\uu}_h - \uu_h)\cdot\ww_h \dif s,
\end{equation}
where $\hat{\uu}_h = 0.5 (\uu_h - \uu^\text{oth}_h)$ is the centered flux.
In general we shall use the upwind flux, and will therefore only specify
if the centered flux is used. We note that the upwind choice is dissipative,
whereas a central flux would not add any dissipation, i.e., is conservative.

Combining the above, the formulation is
\begin{equation*}
\begin{split}
a'(\sigma_h, \tau_h) + b'(\tau_h,\uu_h)&= 0  \quad \forall \tau_h \in \Sigma_h, \\
(\partial_t\uu_h,\vv_h)_{L^2(\Omega)}+ b'(\sigma_h,\vv_h) + c_{upw}(\uu_h,\uu_h,\vv_h)
     + b(\vv_h, p_h) &= 0 \quad \forall \vv_h \in V_h,\\
b(\uu_h, q_h) &= 0 \quad \forall q_h \in Q_h.
\end{split}
\end{equation*}
Here, $\Sigma_h$ is a space of discontinuous piecewise polynomials of degree $k$,
whose values are trace-zero tensors satisfying continuity of the
`normal-tangential component' across edges and where the normal-tangential
component is a polynomial of order $k-1$ on each edge,
$V_h$ denotes the BDMk finite-element space of order $k$,
and $Q_h$ consists of discontinuous piecewise polynomials of degree $k-1$.

As $\hdiv$ spaces are characterized by normal continuity across
facets, the boundary condition for the normal components are imposed strongly.
The implementation uses hybridisation for the tangential component of the
velocity. That is, a second finite element space is used to couple the
tangential component across element facets and the tangential part of
the boundary condition is enforced strongly on the facet space.

For time-stepping we use a second-order implicit-explicit (IMEX) scheme known as SBDF2~\cite{ARW95}.
This scheme uses the BDF2 scheme to discretise the time-derivative, the Stokes
term is fully implicit, and the non-linear convective term is treated explicitly
with a second order extrapolation. This results in the fully discrete problem to find  $(\sigma_h,\uu_h,p_h)\in\Sigma_h\times V_h\times Q_h$, such that
\begin{equation*}
\begin{multlined}
   a'(\sigma_h^n, \tau_h) + b'(\tau_h,\uu_h^n) + b'(\sigma_h^n,\vv_h)
   + {\textstyle \frac{1}{\Delta t}}(\threehalfs\uu_h^n - 2 \uu_h^{n-1} + \half\uu_h^{n-2},\vv_h)_{L^2(\Omega)}\\
    + 2 c_{upw}(\uu_h^{n-1},\uu_h^{n-1},\vv_h) - c_{upw}(\uu_h^{n-2},\uu_h^{n-2},\vv_h)
    + b(\vv_h, p_h^n) +b(\uu_h^n, q_h) = 0,
\end{multlined}
\end{equation*}
holds for all $(\tau_h,\vv_h,q_h)\in\Sigma_h\times V_h\times Q_h$.
To regularize the saddle-point system and to fix the pressure constant in the
case of pure Dirichlet boundary conditions, we may add a small perturbation to
the lower right pressure block of the block-system resulting from
the above equation. That is, we modify the equation for the divergence constraint to
\begin{equation}\label{eqn.perturbed-div-constrant}
   -\int_\Omega q_h (\nabla \cdot \uu_h^n) + \epsilon \int_\Omega p_h^n q_v\dif\xx
   = -\int_\Omega q_h ((\nabla \cdot \uu_h^n) - \epsilon p_h^n)\dif\xx
   =0,
\end{equation}
with some $\epsilon > 0$. Since our discrete spaces $V_h, Q_h$ satisfy the
discrete de-Rahm complex~\cite{JLM17}, we have the property that
$\nabla\cdot V_h = Q_h$. Consequently, we can replace $q_h$ with
$\nabla\cdot\vv_h$ in \eqref{eqn.perturbed-div-constrant}. This gives us an
equation for the pressure, which we may plug into the second equation in
the above saddle-point system. This leads to the discrete system,
\begin{equation}\label{eqn.hdiv.reduced-system}
\begin{split}
   a'(\sigma_h^n, \tau_h) + b'(\tau_h,\uu_h^n)
   + {\textstyle \frac{1}{\Delta t}}(\threehalfs\uu_h^n,\vv_h)_{L^2(\Omega)}
   + b'(\sigma_h^n,\vv_h) - b(\vv_h, \textstyle{\frac{1}{\epsilon}\nabla\cdot\uu_h}).
\end{split}
\end{equation}
This is a smaller linear system to solve, and the system is no longer a saddle
point problem, making it an easier system to solve numerically. We want to
choose $\epsilon$ small to keep the perturbation from the divergence constraint
small. However, this means that the $\frac{1}{\epsilon}$-scaled term causes the
linear system \eqref{eqn.hdiv.reduced-system} to be very ill-conditioned.
Consequently, we need to use direct linear solvers to handle these
ill-conditioned linear systems, which in turn limits us with respect to the
maximal size of the system we can solve and the parallel scalability.
With the small perturbation, the divergence of $\uu_h$ is no longer exactly zero,
but $\norm{\sdiv\uu_h^n}_{L^2}\leq\epsilon \norm{p_h^n}_{L^2}$.
In the computations, we chose $\epsilon=10^{-12} / \nu$.

We will refer to the approach above as $\MCSfem_k$, where $k$ is the order of
the velocity finite element space. When we use the central flux in the
convective term, we refer to the method as $\MCSfem_k^\text{cf}$.
These approaches are implemented using
Netgen/NGSolve \cite{schoeberl_netgen,schoeberl_ngsolve} and based on a template
available on github\footnote{\url{https://github.com/NGSolve/modeltemplates}}.

\begin{remark}[Other $\hdiv$-conforming methods]
We note there are other $\hdiv$-conforming mass conserving finite element
methods that do not reformulate the problem into a first-order system.
These include \cite{CKS06,SL17a} and hybridized approaches, e.g., \cite{LS16,LLS18a}.
\end{remark}

\subsection{Taylor--Hood}
To compare with more established methods, we also consider Taylor--Hood elements
consisting of $H^1$-conforming finite elements of order $k\geq 2$ for the
velocity and order $k-1$ for the pressure.

To stabilize this scheme, we optionally add grad-div stabilization to improve
the conservation of mass~\cite{OLSHANSKII20093975, dFGJN17}. That is, we add the
term
\begin{equation}
   j_h(\uu_h,\vv_h) \coloneqq \gamma_{gd}\int_\Omega\nabla\cdot\uu_h \nabla\cdot\vv_h\dif \xx,
\end{equation}
to our bilinear form. In our computations, we take $\gamma_{gd}=10^3$, and we
refer to this method as $\gdTHfem_k$, where $k$ is the order of the velocity space.
In case no grad-div stabilization is used, we refer to the method as $\THfem_k$.
Furthermore, we note that it is known that for $\gamma_{gd}\rightarrow\infty$,
the grad-div stabilized Taylor--Hood solution converges to the Scott--Vogelius
solution~\cite{ref:LinkeRebholzSVconnection,LINKE2011612}.

For the convective term, there are multiple choices possible. We take the
divergence form~\cite{JohnBook16}
\begin{equation}
   c_{div}(\uu,\vv,\ww) = c_{conv}(\uu,\vv,\ww)
      + \frac{1}{2}\int_{\Omega}(\nabla\cdot\uu)\vv\cdot\ww\dif\xx,
\end{equation}
which is equivalent to the standard convective form on the continuous level,
but vanishes under less restrictive conditions for $\vv=\ww$ on the discrete
level.

Thus the variational formulation takes the form
\begin{equation}\label{eqn:thfarvorm}
\begin{split}
   \begin{multlined}[b]
      (\partial_t\uu_h,\ww_h)_{L^2(\Omega)}+ a(\uu_h, \ww_h)
      + b(\ww_h,p_h) + j_h(\uu_h,\ww_h)\\
      + c_{div}(\uu_h,\uu_h,\ww_h)
   \end{multlined}
   &= 0 \quad \forall \ww_h \in V_h, \\
b(\uu_h, q_h) &= 0 \quad \forall q_h \in Q_h,
\end{split}
\end{equation}
where $V_h$ is the standard Lagrange space of continuous,
piecewise polynomials of degree $k\geq 2$, and
$Q_h$ is the standard Lagrange space of continuous,
piecewise polynomials of degree $k-1$.

To discretize in time-we use the fully implicit Crank--Nicolson scheme, that is we solve for $(\uu_h,p_h)\in V_h\times Q_h$, such that
\begin{equation}\label{eqn:thfarvorm-cn}
\begin{multlined}
   {\textstyle\frac{1}{\Delta t}}(\uu_h^n-\uu_h^{n-1},\ww_h)_{L^2(\Omega)}
   +\half\big(a(\uu^n_h,\ww_h) + c_{div}(\uu^n_h,\uu^n_h,\ww_h)\big)
   + b(\ww_h,p^n_h) \\
   + j_h(\uu^n_h,\ww_h)
   + \half\big(a(\uu^{n-1}_h, \ww_h) +c_{div}(\uu^{n-1}_h,\uu^{n-1}_h,\ww_h)\big)
   + b(\uu^n_h, q_h) = 0,
\end{multlined}
\end{equation}
for all $(\vv_h,q_h)\in V_h\times Q_h$.
Note that we do not take into account the pressure from the previous time step,
nor do we add the grad-div stabilization term from the last time step. This is
because the method needs to correctly reflect that the pressure is the
Langrange-multiplier to weakly enforce the divergence constraint \cite[Section 4.1.4]{Ric17}.
Similarly, the role of the grad-div form is to penalize deviations from the
divergence constraint.

We use a quasi-Newton method with line search to solve the nonlinear systems
arising in every time step. More precisely, we only recompute the factorization
of the inverse Jacobian if the residual decreases by less than a
factor of $0.3$. The line search consists of a check that the Newton update
has decreased the residual. Otherwise, we reduce the size of the update by half
recursively until the residual has decreased compared to the previous iteration.
We start each quasi-Newton iteration with $(u_h,p_h)$ from the previous time step.
As the scheme is fully implicit,
we can use time steps independent of the mesh size, and as the Jacobian
rarely has to be refactorized, the method is essentially as efficient as
IMEX schemes.

The Taylor--Hood methods are also implemented using Netgen/NGSolve. The meshes
used are constructed as those for the $\hdiv$-conforming method in
Section~\ref{sec:disc.subsec:hdiv}.

\subsection{Scott--Vogelius}

As a final method, we also consider the Scott--Vogeluis finite element pair.
This finite element pair is very similar to the Taylor--Hood pair, using
$H^1$-conforming finite elements for the velocity but discontinuous elements
of order $k-1$ for the pressure. Consequently, we use the same variational
formulation time-stepping scheme used for the unstabilized Taylor--Hood method.
For $k\geq 4$ in two dimensions, this method is inf-sup stable on meshes
without nearly singular vertices~\cite{Vog83,SV85} and for $k\geq d=2$ on
barycentric refined meshes~\cite{Qin94,Zha04}. This method is also exactly
divergence-free and pressure-robust.

There are different ways that the Scott--Vogelius method \cite{lrsBIBih}
can be implemented.
In \cite{lrsBIBki}, the iterated penalty method was used to avoid explicit
knowledge of the pressure space. However, this approach was plagued by nearly
singular vertices \cite{lrsBIBkj} in the meshes used.
Here, we utilize the Netgen mesh generator, which, for the cylinder problem,
does not have any singular or nearly singular vertices.
Thus, the mixed-method implementation of Scott--Vogelius can be used, and in
this case, there is only one change from the (unstabilized) Taylor--Hood method,
namely, the pressure elements are discontinuous instead of continuous.
One drawback of this approach is the linear systems are larger, and it
appears that for Reynolds number 8000 (\Cref{tabl:reeighthousand}),
it is no longer possible to keep the divergence very small.
Upon closer examination, it was found that the Newton step was not converging
at a few time steps, and the divergence also spiked there.
In general, the divergence was very small except at these few times.
Furthermore, it is worth mentioning that, unlike the unstabilized Taylor--Hood
method, the lack of convergence of Newton's method did not lead to a blow-up
of the solution.

\subsection{Meshes}

\begin{figure*}
   \centering
   \includegraphics[height=3cm, trim={0 0 10cm 0}, clip]{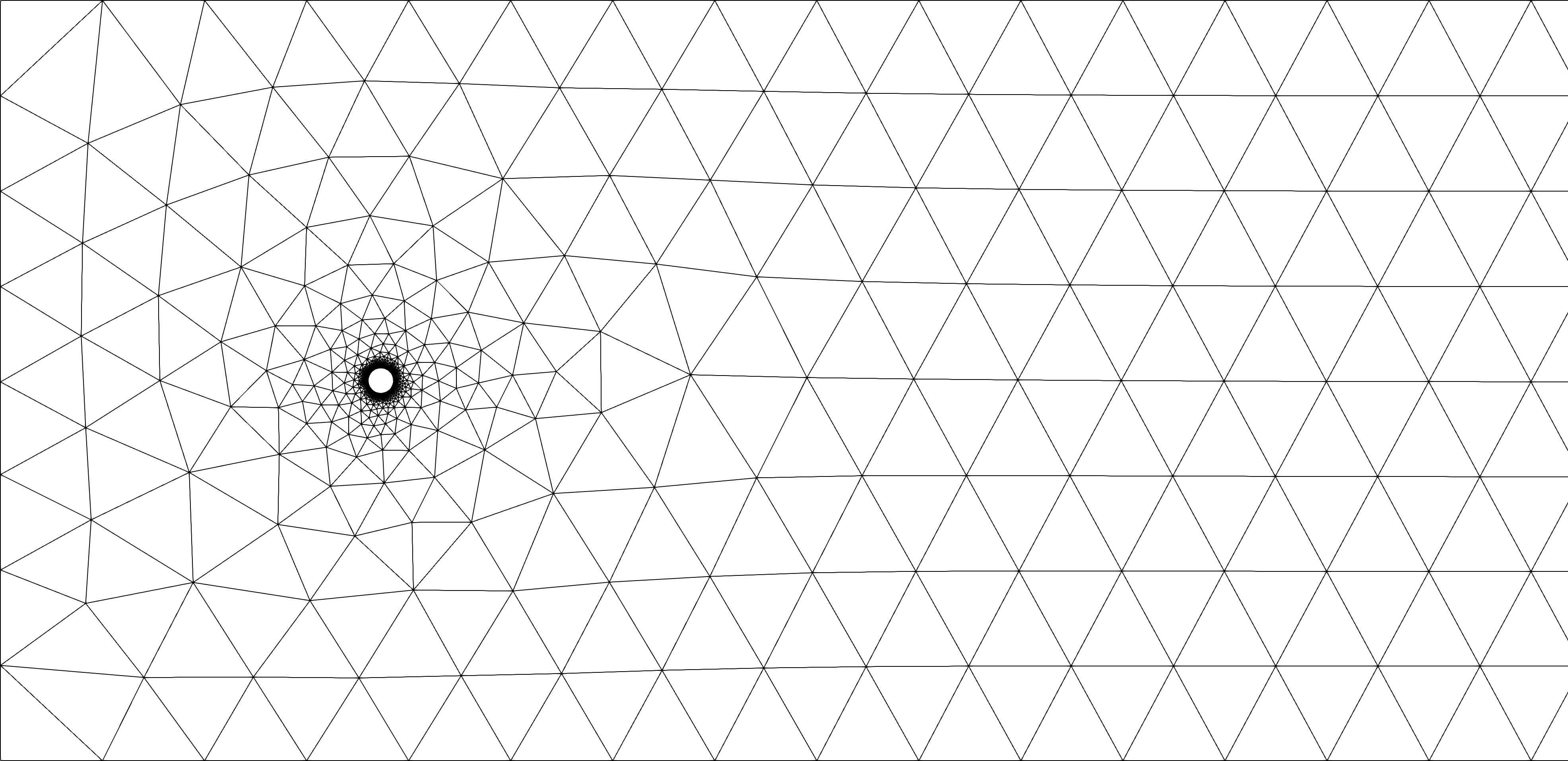}\qquad
   \includegraphics[height=3cm]{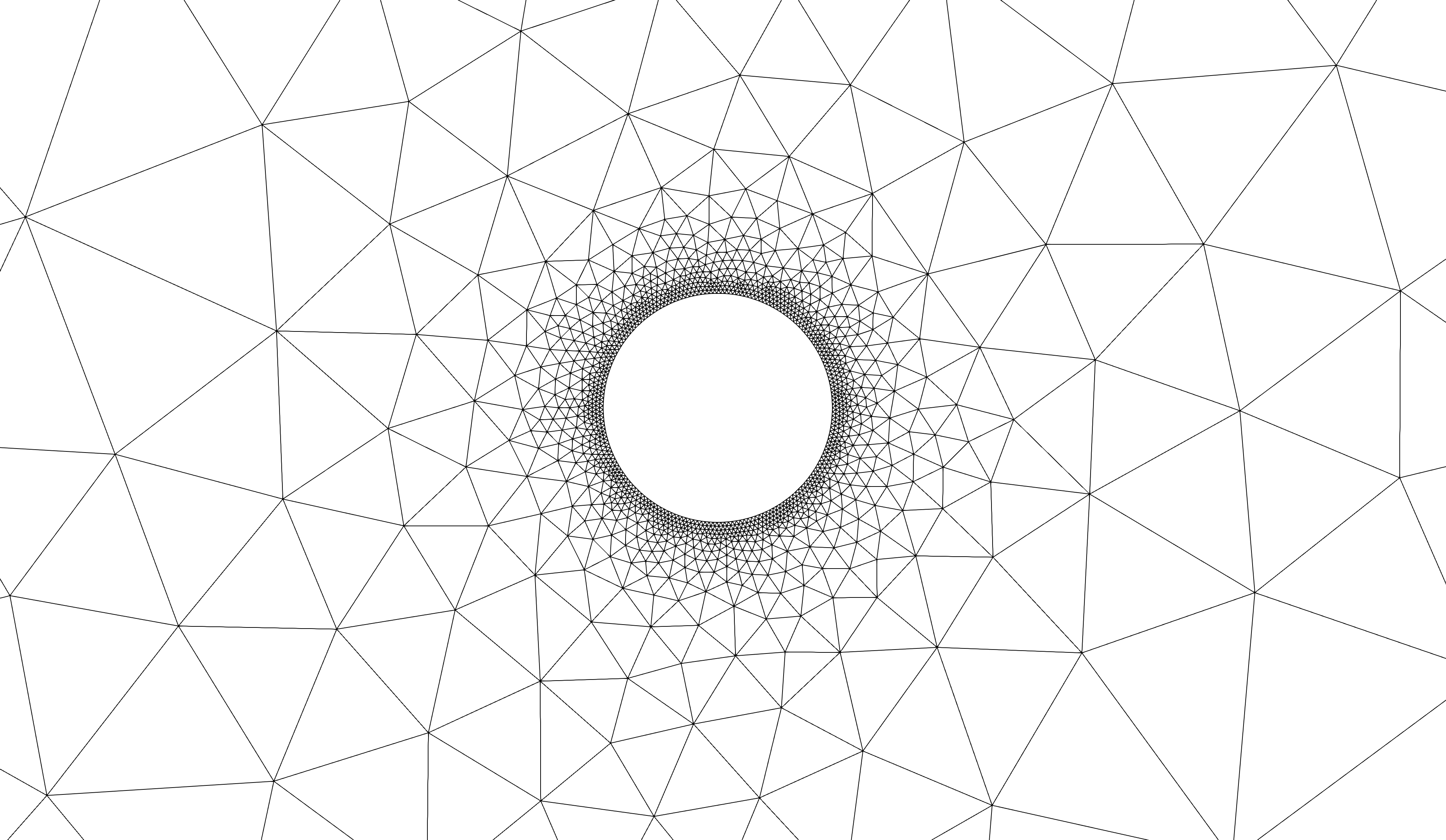}
   \caption{Left: Left section of mesh for the $\THfem{}$, $\SVfem$, and $\MCSfem{}$
      methods with $h=8$. Right: Zoom in to the cylinder of the same mesh.}
   \label{fig:netgen-mesh}
\end{figure*}

We used several meshes for our experiments. For the methods described in
Section~\ref{sec.disctete-methods}, the meshes are unstructured simplicial
meshes constructed using Netgen~\cite{schoeberl_netgen} using a bulk mesh
parameter $h_\text{max}$ and a local mesh size $h_\text{min}=h_\text{max}/250$
on the boundary of the cylinder. The part of the mesh containing the cylinder with
$h_\text{max}=8$ is visualized in \Cref{fig:netgen-mesh}. This is the coarsest
mesh considered.

For $h_\text{max}=8,4,2$, the Taylor--Hood finite element space with $k=4$ has a total of 6994, 163107, and 437439 degrees of freedom (dofs), respectively, for the
saddle point system. The Scott--Vogelius saddle point systems have 87858, 205176,
and 551434 dofs, respectively, due to the larger pressure space. On the same
meshes, the MSC method with $k=4$ has 95856, 224487, and 604047 dofs for the stress
and velocity space in which we solve the problem. However, we emphasize that a
direct comparison of dofs is not useful to compare methods, as different types
of systems must be solved for each method, and the methods have significant
differences regarding their approximation power.

\subsection{Computing the period}

Following \cite{lrsBIBki}, we estimate the period of the drag-lift trajectory
as follows. For selected points along the trajectory (about 100 per period), we
consider points from roughly one-half period away to three-halfs away and
determine the point of closest approach. We start with an estimated average
period, and for each initial point we, obtain an estimated period. Histograms
of such estimates are given in \Cref{fig:histoper}. Taking the average
period from this data gives an estimate of the Strouhal period that can be used
to refine the process, providing a new estimated average period. To allow for
estimates of the period at a finer grain than the individual time steps, we
used a quadratic interpolation of estimated periods to allow for prediction of
the period to second order. That is, we found the trajectory error at the
optimal point and then considered the error at the time steps on either side of
this. Fitting this data with a quadratic and finding the minimum of the
quadratic provided an estimated period to within a fraction of $\Delta t$.

\section{Computational experiments}
\label{sec.numex}

We consider Reynolds numbers between $120$ and $8000$. The results for Reynolds
numbers 120, 250, 500, 1000, 2000, 4000, and 8000 are shown in
\Cref{tabl:reonetwenty,tabl:retwofifty,tabl:refivehundred,tabl:rethousand,%
tabl:retwothousand,tabl:refourthousand,tabl:reeighthousand}. Furthermore, results
for Reynolds numbers between 1100 and 1800 are shown in \Cref{tabl:nmvesitived}.
Finally, the vorticity of the resulting flow computed on the finest mesh
with the MSC method is shown for Reynolds numbers 500, 1000, and 2000
in \Cref{fig.vorticity.Re500,fig.vorticity.Re1000,fig.vorticity.Re2000}.

\subsection{Reynolds 120--500}

We began with a small Reynolds number (120) to provide a low bar for diverse
methods. Nevertheless, \Cref{tabl:reonetwenty} shows that unstabilized,
low-order ($k=2$) Taylor--Hood, which is commonly thought of as the workhorse
for fluid flow, struggled to get a viable result. By contrast, all of the
high-order methods agree to several digits in the prediction of the Strouhal
period and drag. The stabilized Taylor--Hood for $k=2$ is between the
high-order predictions and the unstabilized $k=2$ results. Note that at this
Reynolds number, stabilization of high-order ($k=4$) Taylor--Hood does not make
a difference, except on the size of the divergence error.

\Cref{tabl:retwofifty} tells a similar story with Reynolds number 250.
The high-order methods are in agreement and the low-order methods are
substandard, in particular suggesting erroneous chaotic behavior. This chaotic behavior was also observed in \cite{lrsBIBki}, leading to the suggestion that the flow should become chaotic at lower Reynolds numbers compared to some of the literature. We now understand this to be due to under-resolved simulations in \cite{lrsBIBki}, where high-frequency errors get amplified by the non-linear term.

\begin{figure}
   \centering
   \includegraphics[width=\textwidth]{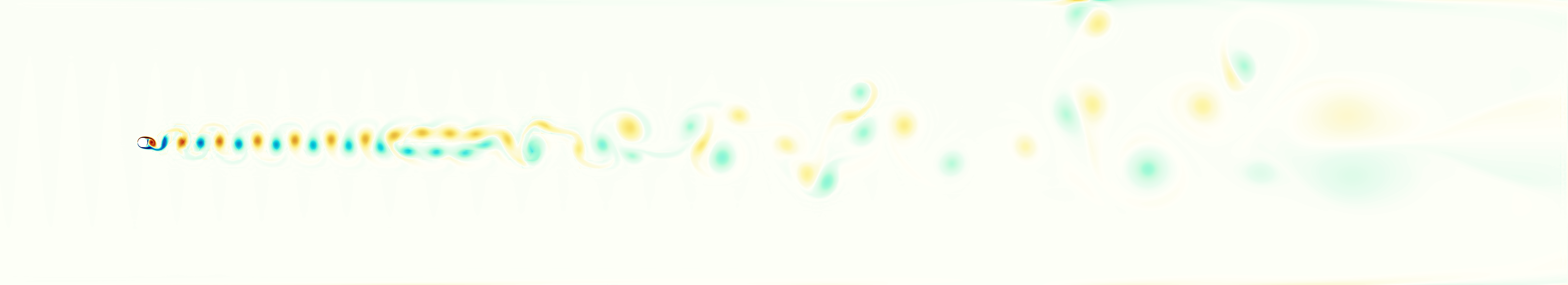}
   \\[3pt]
   \includegraphics[height=.8cm]{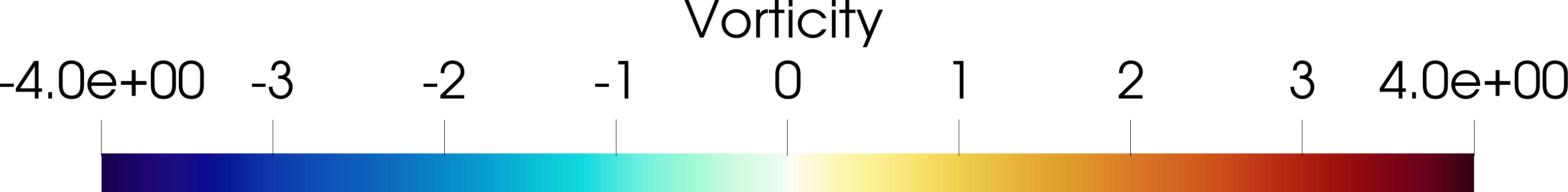}
   \caption{Vorticity of the flow with Reynolds number 500 at $t=300$ computed
      using the $\MCSfem_{4}$ method on the mesh with $h_\text{max}=2$.}
   \label{fig.vorticity.Re500}
\end{figure}

For Reynolds numbers 500 and higher, we decided to drop the lowest-order
Taylor--Hood simulation data. This is due to the method's poor performance
at the lower Reynolds numbers. In some cases, the unstabilized method did not
even continue to time 500, with the Newton iteration failing and the solution
blowing up. For Reynolds number 500, c.f.~\Cref{tabl:refivehundred}, even
the stabilized lowest-order Taylor--Hood simulation data is not viable.
For higher Reynolds numbers, we no longer report results using this method.
Finally, the vorticity of the velocity solution from the MCS method on the
finest mesh is shown in \Cref{fig.vorticity.Re500}. Here, we see the
periodic flow in the wake of the cylinder for about 70 spatial units, after
which the flow transitions away into more chaotic behavior. However, since our
mesh is very coarse in this part of the domain compared to the resolution
around the cylinder, c.f.~\Cref{fig:netgen-mesh}, we do not claim that this is the fully resolved behavior
so far downstream in the channel.

\begin{table}
   \centering
   \caption{Results for Reynolds number 120. The consensus values for the average drag is 1.348
      and for the Strouhal period is 11.34.}
   \label{tabl:reonetwenty}
   \begin{tabular}{rcllllll}
      \toprule
      Method & $h_\text{max}/h_\text{min}$ & $\Delta t$ & $t_{start}$ & $t_{end}$ & Drag & Period & $\Vert\nabla\cdot\uu_h\Vert_{\ell^\infty(L^2)}$\\
      \midrule
      $\MCSfem_4$ & 8 / 0.032 & 0.005   & 280 & 480 & 1.3486 & $11.339\pm 0.00073$ & NA\\
      $\MCSfem_4$ & 4 / 0.016 & 0.0025  & 280 & 480 & 1.3485 & $11.337\pm 0.00085$ & NA\\
      $\MCSfem_4$ & 2 / 0.008 & 0.00125 & 280 & 480 & 1.3482 & $11.337\pm 0.00098$ & NA\\
      \cmidrule(lr){1-8}
      $\SVfem_{4}$ & 8 / 0.032 & 0.01    & 280 & 480 & 1.3542 & $11.303\pm 0.00354$ & $2.83\times10^{-7}$\\
      \cmidrule(lr){1-8}
      $\gdTHfem_4$ & 8 / 0.032 & 0.01    & 280 & 480 & 1.3541 & $11.304\pm 0.00350$ & $3.44\times10^{-3}$\\
      $\gdTHfem_4$ & 4 / 0.016 & 0.01    & 280 & 480 & 1.3492 & $11.339\pm 0.00005$ & $1.52\times10^{-3}$\\
      $\gdTHfem_4$ & 2 / 0.008 & 0.01    & 280 & 480 & 1.3476 & $11.338\pm 0.00010$ & $5.72\times10^{-4}$\\
      \cmidrule(lr){2-8}
      $\gdTHfem_2$ & 8 / 0.032 & 0.01    & 280 & 480 & 1.4633 & $11.139\pm 0.28907$ & $1.35\times10^{-2}$\\
      $\gdTHfem_2$ & 4 / 0.016 & 0.01    & 280 & 480 & 1.3222 & $11.814\pm 0.05813$ & $1.10\times10^{-2}$\\
      \cmidrule(lr){1-8}
      $\THfem_4$   & 8 / 0.032 & 0.01    & 280 & 480 & 1.3443 & $11.379\pm 0.00049$ & $3.74\times10^{0}$\\
      $\THfem_4$   & 4 / 0.016 & 0.01    & 280 & 480 & 1.3492 & $11.339\pm 0.00005$ & $2.57\times10^{0}$\\
      \cmidrule(lr){2-8}
      $\THfem_2$   & 8 / 0.032 & 0.01    & 280 & 380 & 1.1136 & $9.518\pm 3.01488$ & $8.61\times10^{1}$\\
      $\THfem_2$   & 4 / 0.016 & 0.01    & 280 & 480 & 1.0662 & $7.898\pm 2.89044$ & $5.42\times10^{0}$\\
      $\THfem_2$   & 2 / 0.008 & 0.01    & 280 & 480 & 1.2397 & $9.368\pm 2.66038$ & $2.35\times10^{0}$\\
      \bottomrule
   \end{tabular}
\end{table}

\begin{table}
   \centering
   \caption{Results for Reynolds number 250. The consensus values for the average drag is 1.376
         and for the Strouhal period is 9.680.}
   \label{tabl:retwofifty}
   \begin{tabular}{rcllllll}
      \toprule
      Method & $h_\text{max}/h_\text{min}$ & $\Delta t$ & $t_{start}$ & $t_{end}$ & Drag & Period & $\Vert\nabla\cdot\uu_h\Vert_{\ell^\infty(L^2)}$\\
      \midrule
      $\MCSfem_4$ & 8 / 0.032 & 0.005   & 280 & 480 & 1.3750 & $9.681\pm 0.00011$ & NA\\
      $\MCSfem_4$ & 4 / 0.016 & 0.0025  & 280 & 480 & 1.3769 & $9.673\pm 0.00013$ & NA\\
      $\MCSfem_4$ & 2 / 0.008 & 0.00125 & 280 & 480 & 1.3760 & $9.680\pm 0.00023$ & NA\\
      \cmidrule(lr){1-8}
      $\gdTHfem_4$ & 8 / 0.032 & 0.01    & 280 & 480 & 1.3867 & $9.658\pm 0.01316$ & $5.77\times10^{-3}$\\
      $\gdTHfem_4$ & 4 / 0.016 & 0.01    & 280 & 480 & 1.3749 & $9.679\pm 0.00007$ & $3.84\times10^{-3}$\\
      $\gdTHfem_4$ & 2 / 0.008 & 0.01    & 280 & 480 & 1.3765 & $9.684\pm 0.00009$ & $1.21\times10^{-3}$\\
      \cmidrule(lr){2-8}
      $\gdTHfem_2$ & 8 / 0.032 & 0.01    & 280 & 480 & 1.3531 & $9.460\pm 1.79078$ & $1.67\times10^{-2}$\\
      $\gdTHfem_2$ & 4 / 0.016 & 0.01    & 280 & 480 & 1.3529 & $9.732\pm 0.35677$ & $1.47\times10^{-2}$\\
      \cmidrule(lr){1-8}
      $\THfem_2$   & 2 / 0.008 & 0.01    & 280 & 480 & 1.0144 & $9.488\pm 2.64040$ & $1.22\times10^{2}$\\
      \bottomrule
   \end{tabular}
\end{table}

\begin{table}
   \centering
   \caption{Results for Reynolds number 500. The consensus values for the average drag is 1.448
      and for the Strouhal period is 8.82.}
   \label{tabl:refivehundred}
   \begin{tabular}{rcllllll}
      \toprule
      Method & $h_\text{max}/h_\text{min}$ & $\Delta t$ & $t_{start}$ & $t_{end}$ & Drag & Period & $\Vert\nabla\cdot\uu_h\Vert_{\ell^\infty(L^2)}$\\
      \midrule
      $\MCSfem_4$ & 8 / 0.032 & 0.002  & 280 & 480 & 1.4464 & $8.836\pm 0.00004$ & NA\\
      $\MCSfem_4$ & 4 / 0.016 & 0.001  & 280 & 480 & 1.4499 & $8.811\pm 0.00006$ & NA\\
      $\MCSfem_4$ & 2 / 0.008 & 0.0005 & 280 & 480 & 1.4477 & $8.822\pm 0.00016$ & NA\\
      \cmidrule(lr){1-8}
      $\gdTHfem_4$ & 8 / 0.032 & 0.01   & 280 & 450 & 1.4537 & $8.823\pm 0.43717$ & $7.82\times10^{-3}$ \\
      $\gdTHfem_4$ & 4 / 0.016 & 0.01   & 280 & 480 & 1.4488 & $8.877\pm 0.00796$ & $4.77\times10^{-3}$ \\
      $\gdTHfem_4$ & 2 / 0.008 & 0.01   & 280 & 480 & 1.4483 & $8.811\pm 0.00007$ & $2.85\times10^{-3}$ \\
      \cmidrule(lr){2-8}
      $\gdTHfem_2$ & 8 / 0.032 & 0.01   & 280 & 480 & 1.3838 & $8.966\pm 1.90009$ & $2.09\times10^{-2}$\\
      $\gdTHfem_2$ & 4 / 0.016 & 0.01   & 280 & 480 & 1.3753 & $8.716\pm 1.20772$ & $1.82\times10^{-2}$\\
      \bottomrule
   \end{tabular}
\end{table}

\subsection{Reynolds 1000 and higher}

\begin{figure}
   \centering
   \includegraphics[width=\textwidth]{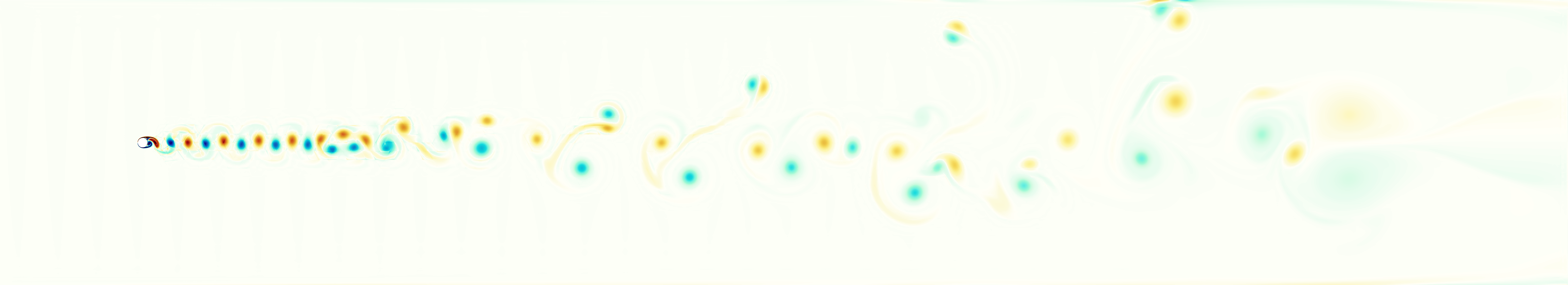}
   \\[3pt]
   \includegraphics[height=.8cm]{img/flow_cylinder_vorticity_colorbar.png}
   \caption{Vorticity of the flow with Reynolds number 1000 at $t=300$ computed
      using the $\MCSfem_{4}$ method on the mesh with $h_\text{max}=2$.}
   \label{fig.vorticity.Re1000}
\end{figure}

For Reynolds numbers 1000 and higher, only the high-order grad-div stabilized
Taylor--Hood, Scott--Vogelius and the MCS method remain viable. To better
understand the behavior of the methods regarding their numerical dissipation,
we also consider the MCS method with the centered flux for the convective term.
This remains an exactly divergence-free but does not include any numeric
dissipation, similar to Scott--Vogelius with Crank--Nicolson time-stepping.
This will allow us to judge the effect of the numerical dissipation from the
upwind convective term in the MCS method used so far, and better judge the
effect of the pressure-robust methods. The results are presented in
\Cref{tabl:rethousand} for Reynolds number 1000. Additionally, the vorticity
of the MSC method on the finest mesh at time $t=300$ is shown in
\Cref{fig.vorticity.Re1000}.

At Reynolds number 1000, the drag-lift dynamics are very close to periodic but
have a very small standard deviation in the computed period. In particular, we
see in \Cref{fig.vorticity.Re1000} that the solution remains periodic for
about 40 units downstream of the cylinder before the solution transitions into
chaos. However, we again cannot claim high accuracy of the solution in this
part of the domain, compared to close to the cylinder. Furthermore, the
$\MCSfem_4^\text{cf}$, $\gdTHfem_4$ and $\SVfem_4$ drag-lift dynamics are only
periodic on the finest mesh considered.

In particular, it is interesting to see in \Cref{tabl:rethousand} that the MCS
method with the upwind convective term results in a periodic solution on all
three meshes. This is in stark contrast to the MCS method with the centered
flux in the convective term, Scott--Vogelius, and grad-div stabilized
Taylor--Hood results. In these methods, the standard deviation of the period
is very large on the coarsest mesh, and the periodic solution is obtained only
using the finest mesh. We note that both the MCS with a centered flux and the
Scott--Vogelius method do not contain any numerical dissipation and that the
grad-div stabilized Taylor--Hood has extra numerical dissipation through the
stabilization term. Interestingly, on the second mesh, the MCS method with
centered flux has a smaller standard deviation than the Scott--Vogelius method,
and the grad-div stabilized Taylor--Hood has the largest standard deviation on
the second mesh. Furthermore, we see that the time-step choice is not the
critical factor here, as the $\gdTHfem_{4}$ results are comparable with those
using a smaller time step.

Overall, it is harder to get a close agreement for the average drag and even
harder for the Strouhal period at this and higher Reynolds numbers. This
illustrates the need for very good methods at higher Reynolds numbers.

\begin{figure}
   \centering
   \includegraphics[width=\textwidth]{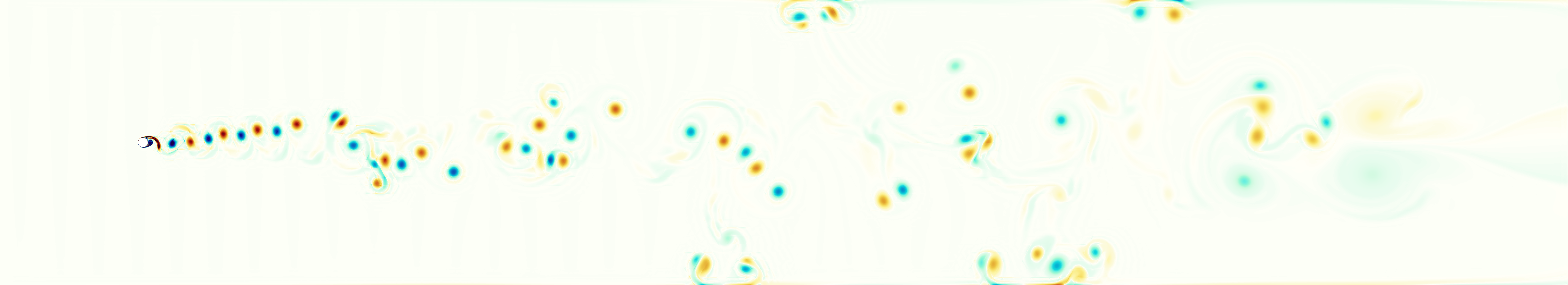}\\[2pt]
   \includegraphics[width=\textwidth]{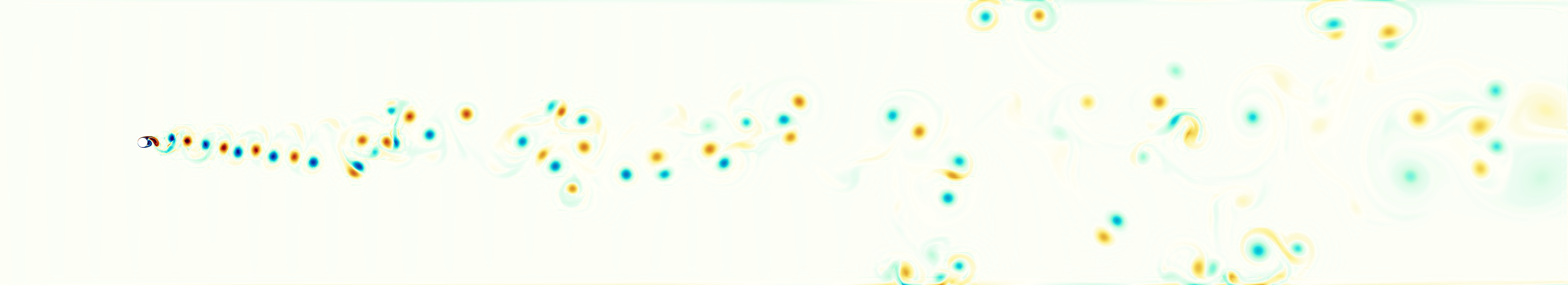}
   \\[3pt]
   \includegraphics[height=.8cm]{img/flow_cylinder_vorticity_colorbar.png}
   \caption{Vorticity of the flow with Reynolds number 2000 computed
      using the $\MCSfem_{4}$ method on the mesh with $h_\text{max}=2$.
      Top: $t=310$, Bottom: $t=350.4$.}
   \label{fig.vorticity.Re2000}
\end{figure}

In \Cref{tabl:retwothousand}, we see for Reynolds number 2000 that the
drag-lift dynamics is no longer periodic, as also observed in
\cite{ref:wolfchaocylinderflow}. Similar
to~\cite[Fig. 7]{ref:wolfchaocylinderflow}, we observe an oscillation of the
wake of the flow behind the cylinder in \Cref{fig.vorticity.Re2000}.

As the drag-lift dynamics is no longer periodic at Reynolds number 2000, we
computed the flow for a sequence of Reynolds number values between 1000 and
2000 to see how the chaos emerges, as shown in \Cref{tabl:nmvesitived}.
We obtain a plausibly periodic solution at Reynolds number 1125, but with an
ever-increasing standard deviation of the Strouhal period as the Reynolds
number increases, even on the finest meshes.

\begin{figure}
   \centering
   \includegraphics{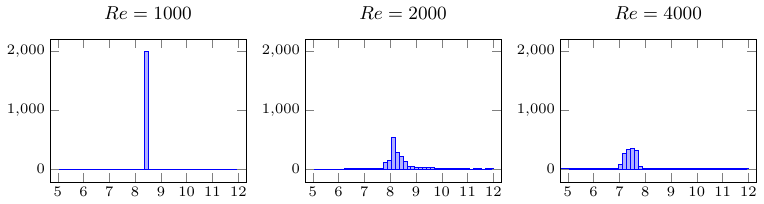}
   \caption{Histogram of the computed periods of the flow at Reynolds numbers
      1000, 2000, 4000 using the $\MCSfem_4$ method with $h_\text{max}=2$.}
\label{fig:histoper}
\end{figure}

For Reynolds numbers 4000 and higher, we no longer attempt to provide consensus
values for the Strouhal period or drag; see \Cref{tabl:refourthousand} and
\Cref{tabl:reeighthousand}. Furthermore, \Cref{fig:histoper}
illustrates how the spread in the computed Strouhal period increases with
increasing Reynolds number.

\begin{table}
   \centering
   \caption{Results for Reynolds number 1000. The consensus values for the average drag is 1.54
      and for the Strouhal period is 8.36.}
   \label{tabl:rethousand}
   \begin{tabular}{rcllllll}
      \toprule
      Method & $h_\text{max}/h_\text{min}$ & $\Delta t$ & $t_{start}$ & $t_{end}$ & Drag & Period & $\Vert\nabla\cdot\uu_h\Vert_{\ell^\infty(L^2)}$\\
      \midrule
      $\MCSfem_4$           & 8 / 0.032 & 0.001   & 280 & 480 & 1.5420 & $8.351\pm 0.00005$ & NA\\
      $\MCSfem_4$           & 4 / 0.016 & 0.0005  & 280 & 480 & 1.5504 & $8.322\pm 0.00036$ & NA\\
      $\MCSfem_4$           & 2 / 0.008 & 0.00025 & 280 & 480 & 1.5384 & $8.361\pm 0.00318$ & NA\\
      \cmidrule(lr){1-8}
      $\MCSfem_4^\text{cf}$ & 8 / 0.032 & 0.0005  & 280 & 480 & 1.5238 & $8.516\pm 1.29090$ & NA\\
      $\MCSfem_4^\text{cf}$ & 4 / 0.016 & 0.0005  & 280 & 480 & 1.5615 & $8.295\pm 0.05556$ & NA\\
      $\MCSfem_4^\text{cf}$ & 2 / 0.008 & 0.00025 & 280 & 480 & 1.5320 & $8.393\pm 0.00879$ & NA\\
      \cmidrule(lr){1-8}
      $\gdTHfem_4$          & 8 / 0.032 & 0.01    & 280 & 480 & 1.5343 & $8.481\pm 1.30527$ & $1.07\times10^{-2}$ \\
      $\gdTHfem_4$          & 4 / 0.016 & 0.01    & 280 & 480 & 1.5162 & $8.297\pm 0.42448$ & $6.79\times10^{-3}$ \\
      $\gdTHfem_4$          & 4 / 0.016 & 0.005   & 280 & 480 & 1.5163 & $8.302\pm 0.41216$ & $6.78\times10^{-3}$ \\
      $\gdTHfem_4$          & 2 / 0.008 & 0.01    & 280 & 480 & 1.5444 & $8.365\pm 0.00169$ & $5.13\times10^{-3}$ \\
      \cmidrule(lr){1-8}
      $\SVfem_4$            & 8 / 0.032 & 0.01    & 280 & 480 & 1.5400 & $8.550\pm 1.31361$ & $9.78\times10^{-6}$\\
      $\SVfem_4$            & 4 / 0.016 & 0.01    & 280 & 480 & 1.5168 & $8.278\pm 0.24277$ & $4.40\times10^{-7}$\\
      $\SVfem_4$            & 2 / 0.008 & 0.01    & 280 & 480 & 1.5444 & $8.365\pm 0.00157$ & $1.18\times10^{-7}$\\
      \bottomrule
   \end{tabular}
\end{table}

\begin{table}
   \centering
   \caption{Results for Reynolds number 2000. The consensus values for the average drag is 1.65
      and for the Strouhal period is 8.4.}
   \label{tabl:retwothousand}
   \begin{tabular}{rcllllll}
      \toprule
      Method & $h_\text{max}/h_\text{min}$ & $\Delta t$ & $t_{start}$ & $t_{end}$ & Drag & Period & $\Vert\nabla\cdot\uu_h\Vert_{\ell^\infty(L^2)}$\\
      \midrule
      $\MCSfem_4$ & 8 / 0.032 & 0.001   & 280 & 480 & 1.6447 & $8.441\pm 0.99875$ & NA\\
      $\MCSfem_4$ & 4 / 0.016 & 0.0005  & 280 & 480 & 1.6470 & $8.417\pm 0.95675$ & NA\\
      $\MCSfem_4$ & 2 / 0.008 & 0.00025 & 280 & 480 & 1.6411 & $8.407\pm 0.94202$ & NA\\
      \cmidrule(lr){1-8}
      $\SVfem_4$   & 2 / 0.008 & 0.01    & 280 & 480 & 1.6603 & $8.066\pm 0.68569$ & $1.75\times 10^{-7}$\\
      \cmidrule(lr){1-8}
      $\gdTHfem_4$ & 4 / 0.016 & 0.01    & 280 & 480 & 1.5798 & $8.902\pm 1.27449$ & $9.55\times 10^{-3}$ \\
      $\gdTHfem_4$ & 2 / 0.008 & 0.01    & 280 & 480 & 1.6524 & $8.365\pm 0.96296$ & $5.61\times 10^{-3}$\\
      \bottomrule
   \end{tabular}
\end{table}

\begin{table}
   \centering
   \caption{Results for Reynolds number 4000.}
   \label{tabl:refourthousand}
   \begin{tabular}{rcllllll}
      \toprule
      Method & $h_\text{max}/h_\text{min}$ & $\Delta t$ & $t_{start}$ & $t_{end}$ & Drag & Period & $\Vert\nabla\cdot\uu_h\Vert_{\ell^\infty(L^2)}$\\
      \midrule
      $\MCSfem_4$ & 8 / 0.032 & 0.001   & 280 & 480 & 1.7778 & $7.867\pm 1.32713$ & NA\\
      $\MCSfem_4$ & 4 / 0.016 & 0.0005  & 280 & 480 & 1.5922 & $9.030\pm 2.33259$ & NA\\
      $\MCSfem_4$ & 2 / 0.008 & 0.00025 & 280 & 480 & 1.7717 & $7.792\pm 1.32658$ & NA\\
      \cmidrule(lr){1-8}
      $\SVfem_4$   & 2 / 0.008 & 0.01    & 280 & 480 & 1.4719 & $8.904\pm 2.68945$ & $1.28\times10^{-6}$\\
      \cmidrule(lr){1-8}
      $\gdTHfem_4$ & 4 / 0.016 & 0.01    & 280 & 480 & 1.6283 & $8.941\pm 2.01980$ & $1.25\times 10^{-2}$\\
      $\gdTHfem_4$ & 2 / 0.008 & 0.01    & 280 & 480 & 1.5620 & $8.997\pm 2.61238$ & $9.48\times 10^{-3}$\\
      \bottomrule
   \end{tabular}
\end{table}

\begin{table}
   \centering
   \caption{Results for Reynolds number 8000}
   \label{tabl:reeighthousand}
   \begin{tabular}{rcllllll}
      \toprule
      Method & $h_\text{max}/h_\text{min}$ & $\Delta t$ & $t_{start}$ & $t_{end}$ & Drag & Period & $\Vert\nabla\cdot\uu_h\Vert_{\ell^\infty(L^2)}$\\
      \midrule
      $\MCSfem_4$ & 4 / 0.016 & 0.0004 & 280 & 480 & 1.6900 & $8.924\pm 2.46617$ & NA\\
      $\MCSfem_4$ & 2 / 0.008 & 0.0002 & 280 & 480 & 1.7262 & $9.252\pm 2.32805$ & NA\\
      \cmidrule(lr){1-8}
      $\SVfem_4$   & 2 / 0.008 & 0.01   & 280 & 480 & 1.6481 & $9.353\pm 2.23443$ & $3.52\times10^{-3}$\\
      \cmidrule(lr){1-8}
      $\gdTHfem_4$ & 4 / 0.016 & 0.01   & 280 & 480 & 1.5174 & $8.883\pm 2.27784$ & $1.29\times10^{-2}$\\
      $\gdTHfem_4$ & 2 / 0.008 & 0.01   & 280 & 480 & 1.6522 & $9.135\pm 2.38565$ & $1.17\times10^{-2}$\\
      \bottomrule
   \end{tabular}
\end{table}

\begin{table}
   \centering
   \caption{Investigating at what Reynolds number the flow becomes chaotic.}
   \label{tabl:nmvesitived}
   \begin{tabular}{crclllll}
      \toprule
      $Re$ & Method & $h_\text{max}/h_\text{min}$ & $\Delta t$ & $t_{start}$ & $t_{end}$ & Drag & Period\\
      \midrule
      1100 & $\MCSfem_4$ & 8 / 0.032 & 0.001   & 280 & 480 & 1.5581 & $8.279\pm 0.00258$ \\
      \cmidrule(lr){2-8}
      1125 & $\MCSfem_4$ & 8 / 0.032 & 0.001   & 280 & 480 & 1.5628 & $8.262\pm 0.00822$ \\
      \cmidrule(lr){2-8}
      1150 & $\MCSfem_4$ & 8 / 0.032 & 0.001   & 280 & 480 & 1.5660 & $8.245\pm 0.02200$ \\
      \cmidrule(lr){2-8}
      1175 & $\MCSfem_4$ & 8 / 0.032 & 0.001   & 280 & 480 & 1.5718 & $8.230\pm 0.07870$ \\
      \cmidrule(lr){2-8}
      1200 & $\MCSfem_4$ & 8 / 0.032 & 0.001   & 280 & 480 & 1.5743 & $8.206\pm 0.33194$ \\
      1200 & $\MCSfem_4$ & 4 / 0.016 & 0.0005  & 280 & 480 & 1.5811 & $8.245\pm 0.91742$ \\
      1200 & $\MCSfem_4$ & 2 / 0.008 & 0.00025 & 280 & 480 & 1.5728 & $8.227\pm 0.48112$ \\
      \cmidrule(lr){2-8}
      1400 & $\MCSfem_4$ & 8 / 0.032 & 0.001   & 280 & 480 & 1.5961 & $8.534\pm 1.07792$ \\
      1400 & $\MCSfem_4$ & 4 / 0.016 & 0.0005  & 280 & 480 & 1.5947 & $8.576\pm 1.17649$ \\
      1400 & $\MCSfem_4$ & 2 / 0.008 & 0.00025 & 280 & 480 & 1.5959 & $8.536\pm 1.11533$ \\
      \cmidrule(lr){2-8}
      1600 & $\MCSfem_4$ & 8 / 0.032 & 0.001   & 280 & 480 & 1.6103 & $8.520\pm 1.15913$ \\
      1600 & $\MCSfem_4$ & 4 / 0.016 & 0.005   & 280 & 480 & 1.6137 & $8.504\pm 1.16444$ \\
      1600 & $\MCSfem_4$ & 2 / 0.008 & 0.0025  & 280 & 480 & 1.6136 & $8.507\pm 1.05176$ \\
      \cmidrule(lr){2-8}
      1800 & $\MCSfem_4$ & 8 / 0.032 & 0.001   & 280 & 480 & 1.6312 & $8.462\pm 1.00298$ \\
      1800 & $\MCSfem_4$ & 4 / 0.016 & 0.005   & 280 & 480 & 1.6305 & $8.475\pm 1.06059$ \\
      1800 & $\MCSfem_4$ & 2 / 0.008 & 0.0025  & 280 & 480 & 1.6327 & $8.481\pm 0.98248$ \\
      \bottomrule
   \end{tabular}
\end{table}

\subsubsection{Just enough numerical dissipation}

As far as we know, there are no direct comparisons between conforming and
non-conforming pressure-robust methods, so our results give some data on this
issue. A direct comparison between fully discontinuous Galerkin and an
$\hdiv$-conforming method was performed in \cite{FKLLS_IJNMF_2019}.
Our results indicate that the $\hdiv$-conforming MCS method with an upwind
stabilized convective term is the most accurate scheme considered here, in
particular at higher-Reynolds numbers. This is in line with the fact that at high Reynolds numbers, some artificial dissipation is needed unless all scales are resolved~\cite{Bur07}. This is due to the non-linear coupling in the Navier--Stokes equations pushing energy (enstrophy in our present, two-dimensional case) to ever higher frequencies. Consequently, for a computational method for high Reynolds number flows to be globally energy consistent with respect to the resolved scales and give optimal resolution where the solution is smooth, it is argued in \cite{Bur07} that the method should have optimal convergence properties for the numerical or artificial dissipation when the solution is smooth and respect the power spectrum where it is rough.

That this is the case for an $\hdiv$-conforming method with upwind stabilization is indicated by the fact that in the inviscid limit, an $\hdiv$-conforming finite element approximation with an upwind flux leads to an optimal error estimate of order $\mathcal{O}(h^{k+0.5})$ \cite{BBG20}, something which is not obtained with the conservative central flux \cite{GSS16}. Furthermore, in the divergence-free and conforming case, extra stabilization is needed to obtain the optimal order error estimate in the inviscid limit \cite{BBCG23}.

Furthermore, we note that in \cite{FKLLS_IJNMF_2019}, it was also concluded that for an $\hdiv$ method, the extra dissipation from the upwind flux in the convective term is of benefit, for three-dimensional under resolved Navier--Stokes flows. Finally, it is argued in \cite{FKG24} that numerical schemes should contain some numerical dissipation in the convective term to overcome problems caused by discretely energy-conserving schemes, such as the inability to represent physical dissipation in the inviscid limit and the resulting accumulation of energy in small scales.

\subsection{Other methods}
\label{sec.numex:subsec.other}

We survey briefly other methods that have been used for simulating flow around
a cylinder.

\subsubsection{Vortex blob}
The publication \cite[Figure 5]{ref:wolfchaocylinderflow} suggests that drag
and lift for flow around a cylinder is nearly periodic at $Re=1000$, in agreement
with our findings, with chaotic flow for $Re\geq 2000$.

\subsubsection{OpenFoam}
OpenFoam may not be the ideal discretization approach for low Reynolds numbers.
In \cite{ref:KornbleuthFoamCylinder}, Strouhal numbers and frequencies
computed with OpenFoam were found to be in error by greater than a factor of two
for $Re\in[55,161]$.

\section{Conclusions}
\label{sec.conclusions}

We considered several finite element methods for solving flow past a cylinder.
All of these methods are well studied both computationally and theoretically.
We found that high-order and pressure-robust methods with the correct amount of
numerical dissipation are necessary to obtain reliable results at higher
Reynolds numbers. The widely-used, lowest-order Taylor--Hood method was not
accurate enough to give appropriate results, even at relatively low Reynolds
numbers. We were able to confirm some earlier results by other researchers.
We found that exact incompressibility and appropriate numerical dissipation are
important at higher Reynolds numbers. This is particularly the case when not
all scales are resolved, which is often the case at high Reynolds numbers. In
particular, we found that both $H^1$-conforming and $\hdiv$-conforming
pressure-robust methods without any numerical dissipation from the convective
term did not lead to accurate results unless very fine meshes where considered.
Finally, we studied in detail how chaos appears to emerge for Reynolds numbers
just above 1000 using the pressure-robust MSC method with upwinding which
appears to have the correct amount of numerical dissipation for reliable results.

\section*{Data Availability Statement}

The code used to realize the presented results, the raw data generated by
this code and videos of key simulations are available on zenodo at
\url{https://doi.org/10.5281/zenodo.11490084}.



\end{document}